\documentclass[11pt]{article}
\usepackage{amssymb,amsfonts,amsmath,amsthm,cite}
\usepackage{graphicx,epsfig,subfigure}
\usepackage{enumerate}
\usepackage{color}

\newtheorem{thm}{Theorem}[section]
\newtheorem{cor}{Corollary}[section]
\newtheorem{lem}{Lemma}[section]

\makeatletter \@addtoreset{equation}{section}

\def\pf{\noindent {\it Proof.\ }}
\def\qed{\hfill \rule{4pt}{7pt}}

\def\red{}

\title{\bf Decycling Number of Linear Graphs of Trees}
\author{Jian Wang$^a$, Xirong Xu$^{b,}$\footnote{\ Corresponding author:
xirongxu@dlut.edu.cn}}
\date{\small
{\it $^a$Department of Mathematics\\
\vskip4pt Taiyuan University of Technology, Taiyuan, 030024, P.R.China\\\vskip5pt
$^b$School of Computer Science and Technology\\ \vskip4pt
Dalian University of Technology, Dalian, 116024, P.R.China}
}

\begin{document}

\maketitle

\noindent {\bf Abstract.} The decycling number of a graph $G$ is the minimum number of vertices whose removal from
$G$ results in an acyclic subgraph. \red{ It is known that determining the decycling number of a graph $G$ is equivalent to finding the maximum induced forests of $G$.} The line graphs of trees are the claw-free block graphs. These graphs have been used by Erd\H{o}s, Saks and S\'{o}s to construct graphs with a given number of edges and vertices whose maximum induced tree is very small. In this paper, we give bounds on the decycling number of line graphs of trees and construct extremal trees to show that these bounds are the best possible. We also give bounds on the decycling number of line graph of $k$-ary trees and determine the exact the decycling number of line graphs of perfect $k$-ary trees.


\allowdisplaybreaks

\section{Introduction}
\label{sec-intro}

Let $G=(V,E)$ be a simple graph, with vertex set $V$ and edge set
$E$. A subset $F\subset V(G)$ is called {\it a decycling set} if the subgraph $G-F$ is acyclic. The minimum cardinality of a decyling  set is called {\it the decycling number} (or {\it feedback number}) of $G$ proposed first by Beineke and Vandell~\cite{bv97}. We use the notation $\nabla(G)$ to denote the decycling number of $G$.

In fact, the problem of determining the decycling
number of a graph is $NP$-complete by Karp~\cite{k72} (also see
\cite{gj79}). The best known approximation algorithm for this
problem has approximation ratio $2$~\cite{bbf99}. Determining the decycling number is difficult even for some elementary graphs. We refer the reader to an original
research paper~\cite{bv97} for some results. Bounds on the
decycling numbers have been established for some
well-known graphs, such as hypercubes~\cite{flp00}, star
graphs~\cite{wwc04}, generalized petersen graphs~\cite{gao15}, distance graphs
and circulant graphs~\cite{8}.

For a graph $G$, let $f(G)$ be the maximum number of vertices in an induced subgraph of $G$ that is a forest. An induced forest with maximum number of vertices is called a {maximum induced forest} of $G$. Determining the decycling number of a graph $G$ is equivalent to finding the maximum induced forest of $G$, since the sum of the two numbers equals the order of $G$.

One can also study induced trees rather than forests in graphs. Let $t(G)$ be the size of maximum induced trees in $G$. The problem of bounding $t(G)$ in the connected graph $G$ was first studied by Erd\H{o}s, Saks and S\'{o}s\cite{ess86} thirty years ago. In their paper, Erd\H{o}s, Saks and S\'{o}s studied the relationship between $t(G)$ and several natural parameters of the graph $G$. They were able to obtain asymptotically tight bounds on $t(G)$ when either the number of edges or the independent number of $G$ were known. Their result showed that $t(G)$ can be very small over graphs with $n$ vertices and $m$ edges.
Given a graph $G$, its line graph $L(G)$ is a graph such that each vertex of $L(G)$ represents an edge of $G$ and two vertices of
$L(G)$ are adjacent if and only if their corresponding edges share a common endpoint in $G$. Erd\H{o}s, Saks and S\'{o}s use line graphs of
trees to construct graphs for which $t(G)$ is surprisingly small. Besides, Erd\H{o}s, Saks and S\'{o}s also considered the problem of estimating the size
of maximum induced tree in $K_r$-free graphs. They use line graphs of regular
trees to construct $K_r$-free graphs for which $t(G)$ is small. Recently, Jacob Fox, Po-Shen Loh and Benny Sudakov improved the results on lower bounds of maximum induced trees in $K_r$-free graphs\cite{fox09}.

A {\it linear forest} in a graph $G$ is a vertex disjoint union of simple paths of $G$.  A {\it maximum linear forest} in $G$ is a linear forest
in $G$ with maximum number of edges. The number of edges in maximum linear forests of graph $G$ is denoted by $l(G)$.  Define the hamiltonian completion number of graph $G$, denoted by $hc(G)$, to be the minimum number of edges that need to be added to make $G$ hamiltonian. The  hamiltonian completion problem was introduced in 1970s by Goodman and Hedetniemi\cite{good74,good75}. Goodman and Hedetniemi\cite{good74} prove the following relation between $l(G)$ and $hc(G)$. For any graph $G$ with $n$ vertices, if $hc(G)>0$, then $l(G)+hc(G)=n$; if $hc(G)=0$, then $l(G)=n-1$.

In this paper, we study the decycling number of line graphs of trees. We show that finding maximum induced forests in line graphs is equivalent to finding maximum linear forests in original graphs.  Let $T$ be a tree on $n$ vertices with diameter $d\geq 4$.  We give lower and upper bounds on $\nabla(L(T))$ as follows. If $d$ is even , then
\[
\left\lceil\frac{n-d-1}{d-1}\right\rceil \leq \nabla(L(T))\leq n-d-1.
\]
If $d$ is odd, then
\[
\left\lceil\frac{n-d-2}{d-2}\right \rceil  \leq \nabla(L(T))\leq  n-d-1.
\]
The extremal line graphs that achieve these bounds are also constructed.

A {\it  $k$-ary tree} is a rooted tree where within each level every node has either 0 or $k$ children. A {\it perfect $k$-ary tree} is a  $k$-ary tree in which all leaf nodes are at the same depth. In this paper, we give bounds on decycling number of line graphs of  $k$-ary trees as follows. Let $T$ be a $k$-ary tree on $n$ vertices. Then
\[
\frac{(k-2)n-k+2}{k}\leq \nabla(L(T)) \leq \frac{(k-1)n-2k+1}{k}.
\]
Moreover, we prove that if $T$ is a perfect $k$-ary tree on $n$ vertices with height $h$, then
\[
\nabla(L(T)) = \frac{(k-1)n-k-(-1)^h}{k+1}.
\]

The rest of this paper is organized as follows. In Section 2, we show that finding maximum induced forests in line graphs  is equivalent to finding maximum linear forests in original graphs. In Section 3, we give lower and upper bounds on the decycling numbers of line graphs of trees with given diameter. In Section 4, we give lower and upper bounds on the decycling number of line graphs of $k$-ary trees.

\section{Maximum Induced Forests in Line Graphs}
\label{sec-2-subgraphs}

In this section, we  prove that the maximum induced forests in line graphs correspond to maximum linear forests in original graphs. Denoted $p(G)$ by the length of the longest paths in $G$.

\begin{lem}\label{lem-lf}
A vertex-disjoint path $P$ in $G$ is longest if and only if $L(P)$ is a maximum induced tree in line graph $L(G)$.
A linear forest $F$ in $G$ is maximum if and only if $L(F)$ is a maximum induced forest in line graph $L(G)$. Thus, $p(G)=t(L(G))$ and $l(G)=f(L(G))$.
\end{lem}
\pf
It is known that if line graphs are claw-free, then they contain no induced $K_{1,3}$. So do their induced trees and induced forests. It follows that every induced tree of a line graph is an induced path and every induced forest of a line graph is an induced linear forest.

Moreover, we shall show that the line graph of a vertex-disjoint path in $G$ is an induced path in $L(G)$ and the induced path in $L(G)$ is also a line graph of a vertex-disjoint path in $G$.
If $P=v_1e_1v_2e_2\ldots v_le_lv_{l+1}$ is a vertex-disjoint path in $G$, in which $v_i$'s are vertices
 and $e_j$'s are edges of $G$. Then, we shall show $L(P)=(e_1,e_2,\ldots, e_l)$ is an induced path in $L(G)$. Otherwise, assume that $(e_j,e_k)$ forms an edge in $L(G)$ and $k>j+1$. Then $e_j$ and $e_k$ share a common ending point in $G$.
 We have $\{v_j,v_{j+1}\}\cap\{v_k,v_{k+1}\}\neq \emptyset$, which contradicts with path $P$ is vertex-disjoint.
Conversely, if $H$ is an induced path in line graph $L(G)$. Let $e_1,e_2,\ldots,e_l$ are $l$ consecutive vertices in $H$.
 Clearly, $P = (e_1,e_2,\ldots,e_l)$ is a path in $G$ and $H=L(P)$.

Thus, a vertex-disjoint path $P$  in $G$  is longest if and only if $L(P)$ is a maximum induced tree in $L(G)$ and a linear forest $F$ in $G$  is maximum if and only if $L(F)$ is a maximum induced forest in $L(G)$.\qed

Clearly, linear forests of $G$ have at most  $n-1$ edges. It implies that $f(L(G))\leq n-1$.  Therefore, we
have the following corollary.
\begin{cor}\label{cor-lh}
For any graph $G$ with $n$ vertices and $m$ edges, $\nabla(L(G))\geq m-n+1$.
\end{cor}

\section{The Decycling Number of Line Graphs of Trees}
\label{sec-extremal}

Let $T$ be a tree with $n$ vertices. An {\it inner vertex}  is a vertex of degree at least two. Similarly, an {\it outer vertex}
(or a {\it leaf}) is a vertex of degree one.
Then the vertices of $T$ can be partitioned into the set of leaves $V_{out}$ and the set of inner vertices $V_{in}$.
The cardinality of $V_{out}$ is denoted by $out(T)$. For any inner vertex $v$, let $ex(v)$ be zero if $v$ has at most two neighbors
 of degree less than three; let $ex(v)$ be $k-2$ if $v$ has $k$ neighbors of degree less than three. Then we have the following lemma.

\begin{lem}\label{lem-lh}
For any tree $T$ with $n$ vertices,
\[
\left\lceil \frac {out(T)+\sum_{v\in V_{in}}ex(v)}{2} \right\rceil \leq hc(T)\leq  out(T)-1.
\]
\end{lem}
\pf
Since $T$ is a tree. By adding $hc(T)$ edges on $T$, we get a hamiltonian graph $G$.
Since $G$ has a hamiltonian cycle $C$, each leaf of $T$ is incident to a new edge in $C$.
For any inner vertex $v$, if $ex(v)$ is greater than zero, then at most two edges that incident to $v$ are in $C$.
It means that at least $ex(v)$ neighbors of $v$ have degree less than or equal to two and do not adjacent to $v$ in $C$.  Each of these neighbors has to be incidence to a new edge  in $C$.
Thus, at least $out(T)+\sum_{v\in V_{in}}ex(v)$ vertices are incidence to new edges. Therefore, we have
\[
hc(T)\geq  \left\lceil \frac {out(T)+\sum_{v\in V_{in}}ex(v)}{2} \right\rceil.
\]

\begin{figure}[ht]
\label{fig-contract}
\center
\includegraphics[scale=0.8]{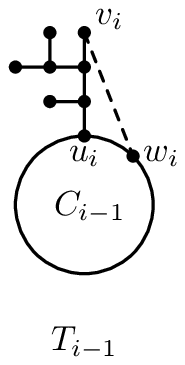}
\hspace{50pt}
\includegraphics[scale=0.8]{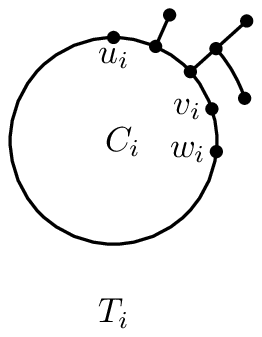}
\caption{An example from tree $T_{i-1}$ to $T_i$.}
\end{figure}

On the other hand, we can get a hamiltonian cycle by adding edges to $T$ according to the following procedure.
Firstly, we choose two leaves $u,v$ of $T$ and add an edge between them. Let $G_0$ be the graph $T+uv$.
Then $G_0$ contains an unique cycle $C_0$ formed by new edge $uv$ and the unique path $P_{u,v}$ in $T$. Let $T_0$ be the graph obtained from $G_0$ by contracting cycle $C_0$, or  $T_0=G_0\cdot C_0$. Denoted by $C_0$ the contracted vertex.
It is easy to see that $T_0$ is a tree with $out(T)-1$. Now choose a leaf $v_1$ outside $C_0$ in $T_0$. Let $P_{C_0v_1}$ be the unique path between $C_0$ and $v_1$. Let $u_1$ be the vertex in the cycle $C_0$ that has the smallest distance to $v_1$. Let $w_1$ be a neighbor of $u_1$ in $C_0$. Then by adding edges $v_1w_1$, we get a larger cycle $C_1=C_0-u_1w_1+P_{u_1v_1}+v_1w_1$. Now let $T_1$ be the graph obtained from $T$ by contracting $C_1$. Then $T_1$ is a tree  with $out(T)-2$ leaves. Now choose a leaf $v_2$ outside $C_1$ from $T_1$. Let $P_{C_1v_2}$ be the unique path between $C_1$ and $v_2$.  Let $u_2$ be the vertex in the cycle $C_1$ that has the smallest distance to $v_2$. Let $w_2$ be a neighbor of $u_2$ in $C_1$.  Then by adding edge $v_2w_2$, we get a larger cycle $C_2=C_2-u_2w_2+P_{u_2v_2}+v_2w_2$ in $T$. Do this procedure repeatedly, through each step we can get a tree $T_i$ from tree $T_{i-1}$ with leaves less than 1(see Fig.1),
the procedure has to be stopped when the contracted tree has only one vertex.
Then we get a hamiltonian cycle by adding $out(T)-1$ edges in $T$. Thus,  $hc(T)\leq  out(T)-1$.
\qed

Since any tree on $n$ vertices have $n-1$ edges. Then $f(L(T)) + \nabla(L(T)) = n-1$. By Lemma \ref{lem-lf}, we know that $l(T)=f(L(T))$. Moreover,  it is true that $l(T)+hc(T)=n$.  Therefore, we have the following corollary.
\begin{cor}\label{cor-lh}
For any tree $T$ on $n$ vertices,
\[
\left\lceil \frac {out(T)+\sum_{v\in V_{in}}ex(v)}{2} \right\rceil -1 \leq \nabla(L(T))\leq  out(T)-2.
\]\end{cor}

Now we introduce an operation on leaves of trees that does not decrease $l(T)$. For any two  leaves $u_i,u_j$ of $T$, suppose their neighbors are $w_i,w_j$.
We define Leaf-Exchange operation on $T$ as removing edge $w_iu_i$ from $T$ and adding edge $u_ju_i$, the obtained tree is denoted by $T[u_i\rightarrow u_j]$.
\begin{lem}\label{lem-operation}
For any two leaves $u_i,u_j$ of $T$,
$l\left(T[u_i\rightarrow u_j]\right)\geq l(T)$.
\end{lem}
\pf Suppose $F$ is a maximum linear forest in $T$. Then
$F-w_iu_i+u_ju_i$ is a linear forest in $T[u_i\rightarrow u_j]$. Thus, we have $l(T)\leq l(T[u_i\rightarrow u_j])$.
\qed

Let $T$ be a tree on $n$ vertices.
The center of a tree is the set of vertices, from which the greatest distance equals to its radius.
Let $v^*$ be one of the center of tree $T$. Then $T$ can be viewed as a rooted tree with root $v^*$. Moreover, we can partition $V(T)$ into sets $V_0(T)$,$V_1(T),\ldots,V_r(T)$, where
$V_i(T)=\{w | d(v^*,w)=i\}$ and $r$ is the radius of $T$. In case of no confusion, $V_i(T)$ is often abbreviated as $V_i$. The vertex in $V_i$  is called  the vertex at depth $i$. Let $V_{\geq 2}=V_2\cup V_3\ldots \cup V_r$. Let $d(T)$ be the diameter of $T$ and $r(T)$ be the radius of $T$. Let $s(T)$ be the number of degree-two vertices in $V_1(T)$. Then, we define three family of rooted trees on $n$ vertices with diameter at most $d$ as follows.
$$
\left\{\begin{array}{ll}
\mathcal{T}_1(n,d)=&
\big\{T \colon |V(T)|=n \mbox{, } d(T)\leq d \mbox{, } r(T)\leq\lceil\frac{d}{2}\rceil\\[5pt]
&\quad\quad\quad\quad \mbox{ and }deg(v)\leq 2 \mbox{ for } v \in V_{\geq 2}\mbox{, }deg(v)\leq 3 \mbox{ for } v \in V_1 \big\},\\[6pt]
\mathcal{T}_2(n,d)=&\{T \colon  s(T)\leq 3, T \in\mathcal{T}_1(n,r)\},\\[6pt]
\mathcal{T}_3(n,d)=&\{T \colon  2\leq s(T)\leq 3, T \in\mathcal{T}_1(n,r)\}.
\end{array}\right.
$$

By the following three lemmas, we shall show that finding the upper bounds for $l(T)$ on all trees is equivalence to finding that on $\mathcal{T}_3(n,r)$.
\begin{lem}\label{t1}
For any tree $T$ on $n$ vertices with diameter $d$, there exists a tree $T'$ in  $\mathcal{T}_1(n,d)$ such that $l(T)\leq l(T')$.
\end{lem}
\pf
Any tree can be viewed as a rooted tree with its center as the root. Suppose to the contrary, there exist trees that we cannot find trees with larger maximum linear forest in $\mathcal{T}_1(n,d)$. Let $T$ be a counterexample with $|V_1|$ maximum. Clearly,  $T$ is not in  $\mathcal{T}_1(n,d)$. Then, $T$ has a vertex
$v$ in $V_1$ such that $deg(v)\geq 4$ or $T$ has a vertex $v$ in $V_{\geq 2}$ such that $deg(v)\geq 3$. We split the proof into two cases as follows.

{\bf Case 1}. $T$ has a vertex $v$ in $V_1$ such that $deg(v)\geq 4$. Assume that $deg(v)=t+1$ and $t\geq 3$. Then $v$
has one neighbor $v^*$ and $t$ neighbors in $V_2$. Let $v_1,v_2,\ldots,v_t$ be these
$t$ neighbors in $V_2$ and $T_1,T_2,\ldots,T_t$ be subtrees of $v$ with root $v_1,v_2,\ldots,v_t$. Let $F$ be a maximum linear forest of $T$.
Then at most two edges of $v^*v,vv_1,vv_2,\ldots,vv_t$ are in $F$. Since $t\geq 3$, there exists one of $vv_1,vv_2,\ldots,vv_t$ that is
not in $F$. Without loss of generality, we assume $vv_t$ is not in $F$. Then by removing edge $vv_t$ from $T$ and adding edge $v^*v_t$, we get a new tree $\bar{T}$ with $d(\bar{T})\leq d(T)$. Clearly, we have $l(T)\leq l(\bar{T})$ since $F$ is also a linear forest of $\bar{T}$. Moreover, $V_1(\bar{T})$ has more vertices than $V_1(T)$.
Since $T$ is the counterexample with $|V_1|$ maximum, we know that $\bar{T}$ is no longer a counterexample. Therefore,
 there exists a tree $T'$ in  $\mathcal{T}_1(n,d)$ such that $l(\bar{T})\leq l(T')$. Then $l(T)\leq l(\bar{T})\leq l(T')$, which
  contradicts with that $T$ is a counterexample.

{\bf Case 2}. $T$ has a vertex $v$ in $V_k(k\geq 2)$ such that $deg(v)\geq 3$. If $deg(v)\geq 4$, then we can get a contradiction by the same argument as
 in Case 1. Thus, we only need to consider the case $deg(v)=3$. Then $v$ has one neighbor $w$ in $V_{k-1}$  and has two neighbors $v_1$ and $v_2$ in
  $V_{k+1}$. $T_1,T_2$ be subtrees of $v$ with root $v_1,v_2$. Let $F$ be a maximum linear forest of $T$. Then at most two edges of $wv,vv_1,vv_2$
  are in $F$.
If $wv$ is not in $F$, by removing edge $wv$ from $T$ and adding edge $v^*v$, we get a new tree $\bar{T}$ with $d(\bar{T})\leq d(T)$. We have $l(T)\leq l(\bar{T})$
 since $F$ is also a linear forest in $\bar{T}$. Since  $V_1(\bar{T})$ is increased by one, $\bar{T}$ is no longer a counterexample.
 Therefore, there exists a tree $T'$ in  $\mathcal{T}_1(n,d)$ such that $l(\bar{T})\leq l(T')$. We get a contradiction. If one of $vv_1$ and $vv_2$
 is not in $F$, without loss of generality, we assume $vv_2$ is not in $F$.
Then by removing edge $vv_2$ from $T$ and adding edge $v^*v_2$, we get a new tree $\bar{T}$, which also leads to a contradiction.

Therefore, the claim holds.
\qed

\begin{lem}\label{t2}
For any tree $T$ in $\mathcal{T}_1(n,d)$, there exists a tree $T'$ in  $\mathcal{T}_2(n,d)$ such that $l(T)\leq l(T')$.
\end{lem}
\pf
Suppose to the contrary, there exist counterexamples.
Let $T$ be the one in $\mathcal{T}_1(n,d)$ with $s(T)$ minimum. Then $v^*$ has at least four degree-two neighbors.
Assume they are $v_1,v_2,\ldots,v_s$. Then by definition of $\mathcal{T}_1(n,d)$, it is easy to see that subtrees with roots $v_1,v_2,\ldots,v_s$ are all paths. Let $P_1,P_2,\ldots,P_s$ be these paths. Clearly, $v_i$ is one endpoint of $P_i$.
 Let $F$ be a maximum linear forest in $T$. Then at least two of edges $v^*v_1,v^*v_2,\ldots,v^*v_s$ are not in $F$.
  Without loss of generality, we suppose $v^*v_1,v^*v_2$ are not in $F$. Then all edges in $P_1,P_2$ are in $F$.

If one of paths $P_1,P_2$ has length at least two.  Without loss of generality, we suppose $P_1$ has length at least 2. Clearly, $v_1$ is one endpoint of $P_1$.
 Let $u_1$ be the other endpoint of $P_1$ and $w_1$ be the parent of $u_1$ in the rooted tree $T$. Then by removing edge $w_1u_1$
  and adding edge $v_1u_1$, we get a new tree $\bar{T}$. We have $l(T)\leq l(\bar{T})$ since $F-w_1u_1+v_1u_1$ is a linear forest of $\bar{T}$.
  Since $s(\bar{T})=s(T)-1$, there exists a tree $T'$ in  $\mathcal{T}_2(n,d)$ such that $l(\bar{T})\leq l(T')$. It leads to a contradiction. If
$P_1,P_2$ all have length one. Assume that $P_1$ is an edge $v_1u_1$ and $P_2$ is an edge $v_2u_2$. Then by removing edge $v_1u_1$ and
 adding edge $v_2u_1$,
we get a new tree $\bar{T}$ with $l(T)\leq l(\bar{T})$, which also leads to a contradiction. Thus, the claim holds.
\qed

\begin{lem}\label{t3}
For any tree $T$ in $\mathcal{T}_2(n,d)$ with $n\geq d$ and $d\geq 4$, there exists a tree $T'$ in  $\mathcal{T}_3(n,d)$ such that $l(T)\leq l(T')$.
\end{lem}
\pf
Suppose $T$ is a counterexample in $\mathcal{T}_2(n,d)$ with $s(T)$ maximum. Clearly, $s(T)\leq 1$. Let $F$ be a maximum linear forest in $T$.

Firstly, we claim $T$ has no leaves in $V_1$.
Otherwise, assume $v_0$ is a leaf in $V_1$. If there is a degree-3 vertex in $V_1$, say $w_0$. Let $u_0$ be a leaf of the subtree with root $w_0$.
Then by Leaf-Exchange operation on $T$, we get a new tree $T[u_0\rightarrow v_0]$ with $s(T[u_0\rightarrow v_0])\geq s(T)+1$ and $l(T[u_0\rightarrow v_0])\geq  l(T)$. Since  $T$ is a counterexample with $s(T)$ maximum. Then there exists a tree $T'$ in $\mathcal{T}_3(n,d)$ such that $l(T[u_0\rightarrow v_0])\leq l(T')$. It follows that  $l(T)\leq l(T')$,  a contradiction. If there is no degree-3 vertex in $V_1$. Then all vertices in $V_1$ have degree at most two. Since $s(T)\leq 1$ and $n\geq d$,  there are at least two leaves in $V_1$, say $v_1$ and $v_2$. Then by Leaf-Exchange operation on $T$, we get a new tree $T[v_1\rightarrow v_2]$ with $s(T[v_1\rightarrow v_2])= s(T)+1$. Since $d\geq 4$, the Leaf-Exchange operation cannot decrease the diameter. Then there exists a tree $T'$ in $\mathcal{T}_3(n,d)$ such that $l(T[v_1\rightarrow v_2])\leq l(T')$. It follows that $l(T)\leq l(T[v_1\rightarrow v_2])\leq l(T')$, a contradiction. Thus, $T$ has no leaves in $V_1$.

Now, we know that there is at most one degree-2 vertex and no leaf in $V_1$. Since $n\geq d\geq 2r(T)$, the number of degree-3 vertices in $V_1$ has to be greater than one.
Let $v_1,v_2,\ldots,v_s$ be degree-3 vertices in $V_1$. If at least one of edges $v^*v_1,v^*v_2,\ldots,v^*v_s$ is in $F$. Without loss of generality,  we assume that $v^*v_1$
is in $F$. Let $u_1,u_2$ be two neighbors of $v_1$ in $V_2$. Then at least one of edges $v_1u_1$ and $v_1u_2$ is not in $F$. Suppose $v_1u_1$ is not in
in $F$. Then by removing edge $v_1u_1$ and adding edge $v^*u_1$, we get a new tree $\bar{T}$ with $l(T)\leq l(\bar{T})$. Clearly, $s(\bar{T})\geq s(T)+1$ and
$\bar{T}$ is in  $\mathcal{T}_2(n,r)$. Then there exists a tree $T'$ in $\mathcal{T}_3(n,d)$ such that $l(\bar{T})\leq l(T')$. It follows that $l(T)\leq l(\bar{T}])\leq l(T')$, a contradiction.  If none of $v^*v_1,v^*v_2,\ldots,v^*v_s$ is in $F$. Then at most one edge incident to $v^*$ is in $F$, and edges in
each subtree with root $v_i$ are all in $F$. Suppose $v_1$ has two neighbors $u_1,u_2$ in $V_2$. Then by removing edge $v_1u_1$ and adding edge $v^*u_1$, we
get a new tree $\bar{T}$. Then $F-v_1u_1+v^*u_1$ is a linear forest in $\bar{T}$. It follows that $l(T)\leq l(\bar{T})$ and $s(\bar{T})\geq s(T)+1$. Then there exists a tree $T'$ in $\mathcal{T}_3(n,d)$ such that $l(\bar{T})\leq l(T')$. It follows that $l(T)\leq l(\bar{T})\leq l(T')$, a contradiction.

Combining all the cases, we complete the proof.
\qed

\begin{thm}\label{thm-ex}
For any tree $T$ on $n$ vertices with diameter $d\geq 4$, we have
$$
\left\{\begin{array}{ll}
 d \leq l(T)\leq  \left \lfloor \frac{(d-2)n+2}{d-1}\right \rfloor,\ for\ d\ is\ even;\\[6pt]
d \leq l(T)\leq   \left \lfloor \frac{(d-3)n+4}{d-2}\right \rfloor,\ for\ d\ is\ odd.
\end{array}\right.
$$
\end{thm}
\pf
For the lower bounds, clearly we have $l(T)=f(L(T))\geq t(L(T))=d(H)$. Moreover, the extremal trees that achieve these lower bounds are shown in Fig.2.

\begin{figure}[ht]
\label{fig-exgraph}
\center
\includegraphics[scale=0.65]{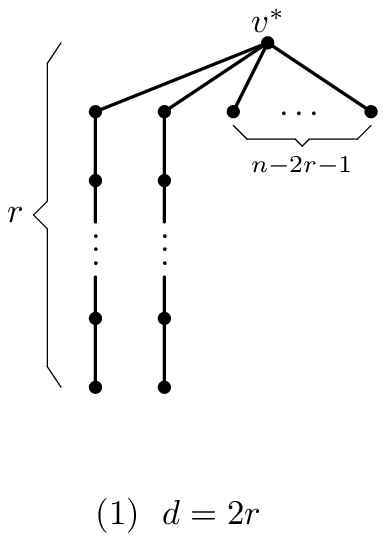}
\hspace{60pt}
\includegraphics[scale=0.65]{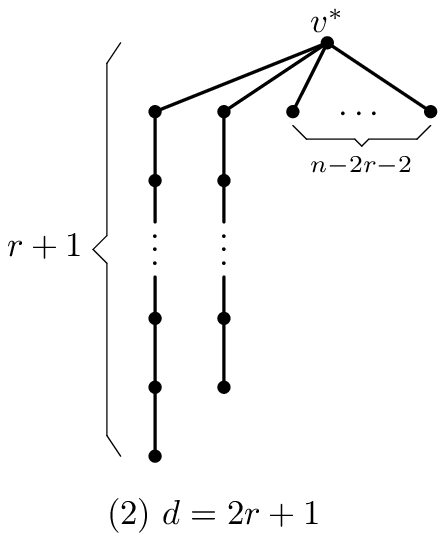}
\caption{Extremal trees that achieve lower bounds.}
\end{figure}

Let $v^*$ be the center of tree $T$. Then $T$ can be viewed as a rooted tree with root $v^*$ and radius $r$. Clearly, $n\geq d$. Then by Lemma \ref{t1}, \ref{t2} and \ref{t3}, we know that there exists a tree $T'\in \mathcal{T}_3(n,d)$ such that $l(T)\leq l(T')$. Therefore,  we only need to consider the upper bounds on $l(T')$ for $T'\in \mathcal{T}_3(n,d)$. Now we split the proof into two cases by the parity of $d$.

\begin{figure}[ht]
\label{fig-exgraph2}
\center
\includegraphics[scale=0.65]{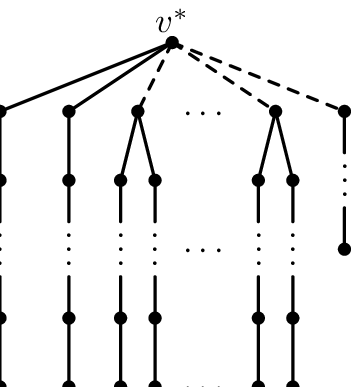}
\hspace{60pt}
\includegraphics[scale=0.65]{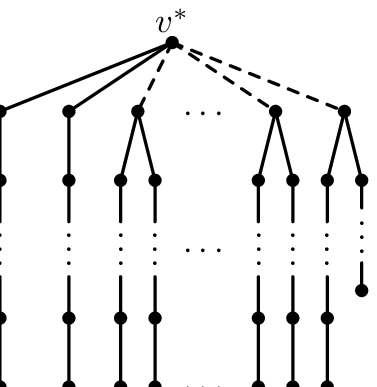}
\caption{Extremal graphs $T^*(n)$ that achieve the upper bounds.}
\end{figure}

{\bf Case 1}. $d=2r$. Let $T$ be a tree in $\mathcal{T}_3(n,d)$.  It is easy to see that $T$ has radius at most $r$. Since $s(T)\geq 2$. Then, at least two vertices in $V_1(T)$ have degree two, say $u$ and $v$. Then the subtree with root $u$ and the subtree with root $v$ are two paths. We call two leaves in these two subtrees critical leaves of $T$ and call all the other leaves non-critical leaves. Let $T^*(n)$ be the tree in $\mathcal{T}_3(n,d)$ satisfying the following two properties as shown in Figure \ref{fig-exgraph2}:

(1) Two critical leaves of $T^*(n)$ are all at depth $r$;

(2) All but at most one of its leaves are at depth $r$. If the only leaf with depth less than $r$ lie in a subtree, whose root is a degree-2 vertices in $V_1(T^*(n))$, then $s(T^*(n))=3$. If the only leaf with depth less than $r$ lie in a subtree, whose root is a degree-3 vertices in $V_1(T^*(n))$, then $s(T^*(n))=2$.

We claim for any tree $T$ in $\mathcal{T}_3(n,d)$, $l(T)\leq l(T^*(n))$. Let $v_1,v_2$ be degree-2 vertices in $V_1$ and $v_3,\ldots,v_k$ be degree-3 vertices in $V_1$. If $s(T)=3$, let $v_{k+1}$ be the third degree-2 vertex in $V_1$. We arrange the subtrees of $T$ with roots $v_1, v_2, v_3, \ldots, v_k, v_{k+1}$ from left to right in the plane. Then do Leaf-Exchange operation from a rightmost leaf to a leftmost leaf with depth less than $r$ convectively. Finally, we shall arrive the tree $T^*(n)$. Since Leaf-Exchange operation from leaf to leaf can never decrease the value of $l(T)$. It follows that $l(T)\leq l(T^*(n))$.

Let $m$ be the remainder of dividing $n-2r-1$ by $2r-1$. We splits the proof into two parts by the value of $m$.

{\bf Case 1.1. } If $1\leq m\leq r$, then $s(T^*(n))=3$ as shown in Fig.3 (1). It is easy to see that there are three vertices in $V_1(T^*(b))$. It follows that $ex(v^*)=1$. Moreover, the number of leaves in $T^*(n)$ can be computed as follows.
 \begin{align*}
 out(T^*(n))&=2\left \lfloor \frac{n-2r-1}{2r-1}\right\rfloor+3\\
 &=2\left \lfloor \frac{n-2}{2r-1}\right\rfloor+1
 \end{align*}
By Lemma \ref{lem-lh}, we have
 \begin{align*}
l(T^*(n))&\leq n-hc(T^*(n)) \\[5pt]
&\leq n- \left\lceil \frac {out(T^*)+\sum_{v\in V_{in}}ex(v)}{2}\right\rceil  \\[5pt]
&  = n- \left \lfloor \frac{n-2}{2r-1}\right\rfloor-1\\[5pt]
& = \left \lceil \frac{(2r-2)n+2}{2r-1}\right\rceil-1= \left \lfloor \frac{(2r-2)n+2}{2r-1}\right\rfloor.
\end{align*}
However, by removing the dashed edges as shown in Fig.3 (1), we get a linear forest with
$n-1-\left \lceil \frac{n-2r-1}{2r-1}\right\rceil= \left \lfloor \frac{(2r-2)n+2}{2r-1}\right\rfloor$ edges.
Thus, $l(T^*(n))= \left \lfloor \frac{(2r-2)n+2}{2r-1}\right\rfloor$.

{\bf Case 1.2. } If $m$ is $0$ or between $d+1$ and $2d-2$, then  $T^*(n)$ is  as shown in Fig.3 (2).
For the second case,
we have $ex(v^*)=0$ and $out(T^*(n))=2\left \lceil \frac{n-2r-1}{2r-1}\right\rceil+2=2\left \lceil \frac{n-2}{2r-1}\right\rceil$. By Lemma \ref{lem-lh},
we have
 \begin{align*}
l(T^*(n))&\leq n-hc(T^*(n))\\[5pt]
&\leq n- \left\lceil \frac {out(T^*(n))+\sum_{v\in V_{in}}ex(v)}{2} \right\rceil\\[5pt]
& =  \left \lfloor \frac{(2r-2)n+2}{2r-1}\right\rfloor.
\end{align*}
However, by removing the dashed edge as shown in Fig.3 (2), we get a linear forest with
$n-1-\left \lceil \frac{n-2r-1}{2r-1}\right\rceil= \left \lfloor \frac{(2r-2)n+2}{2r-1}\right\rfloor$ edges.
Thus, $l(T^*(n))=  \left \lfloor \frac{(2r-2)n+2}{2r-1}\right\rfloor$.

Combining the two subcases, we prove that for any tree $T$ with $n$ vertices and diameter $2r$ $(r\geq 2)$,
\[
 l(T)\leq l(T^*)\leq   \left \lfloor \frac{(2r-2)n+2}{2r-1}\right\rfloor.
\]

\begin{figure}[ht]
\label{fig-exgraph3}
\center
\includegraphics[scale=0.65]{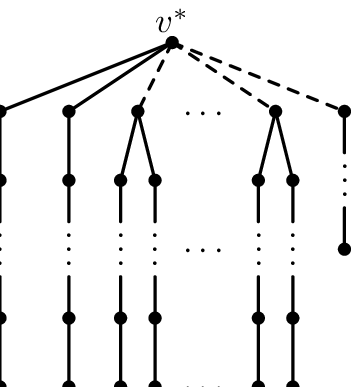}
\hspace{60pt}
\includegraphics[scale=0.65]{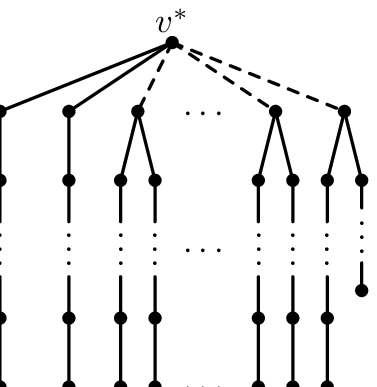}
\caption{Extremal trees $T_1^*(n)$ that achieve the upper bounds.}
\end{figure}

{\bf Case 2}. $d=2r+1$. Let $T$ be a tree in $\mathcal{T}_3(n,d)$. Since $d(T)\leq d$ and $r(T)\leq r+1$. Then if there are leaves in $V_{r+1}$, these leaves are all in the same subtree of $v^*$. Thus, we have $|V_{r+1}|=1$ or $|V_{r+1}|=2$.

{\bf Case 2.1}. $|V_{r+1}|=1$. Let $u$ be the only leaf in $V_{r+1}$ and $w$ be its neighbor. By remove this vertex, we get a tree $T'$ with diameter $2r$.
Suppose $F$ is a maximum linear forest in $T$. Then $F-uw$ is a linear forest in $T'$. Conversely, if $F$ is a maximum linear forest in $T'$ then $F+uw$ is a linear forest in $T'$.  It follows that $l(T)=l(T')+1$. Therefore, by adding one leaf to $T^*(n-1)$, we obtain a new tree $T_1^*(n)$ with $n$ vertices and diameter $2r+1$ as shown in Fig.4. And it is easy to see that for any $T$ in $\mathcal{T}_3(n,d)$, $l(T)\leq l(T_1^*(n))$. Thus,
 \begin{align*}
l(T) &\leq l(T_1^*(n))\\[5pt]
&= l(T^*(n-1))+1\\[5pt]
&= \left \lfloor \frac{(2r-2)n+3}{2r-1}\right\rfloor.
\end{align*}

{\bf Case 2.2}. $|V_{r+1}|=2$. Let $u_1$ and $u_2$ be these two vertices and $w_1$ and $w_2$ be their neighbors, respectively. Clearly, $n\geq 4r+2$.
By remove these two vertices, we get a tree $T'$ with diameter $2r$.
Clearly, we have $l(T)=l(T')+2$. Thus, by adding two leaves to $T^*(n-1)$, we obtain a tree $T_2^*(n)$ with $n$ vertices and diameter $2r+1$ as shown in Fig.5. And any tree T in $\mathcal{T}_3(n,d)$ with
 $|V_{r+1}|=2$ has $l(T)\leq l(T_2^*(n))$. Thus,
 \begin{align*}
l(T)&\leq l(T_2^*(n)) \\[5pt]
&= l(T^*(n-2))+2\\[5pt]
&= \left \lfloor \frac{(2r-2)n+4}{2r-1}\right\rfloor.
\end{align*}

\begin{figure}[ht]
\label{fig-exgraph4}
\center
\includegraphics[scale=0.65]{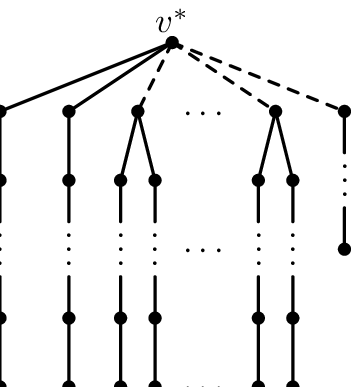}
\hspace{60pt}
\includegraphics[scale=0.65]{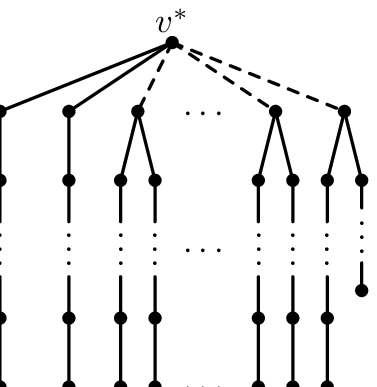}
\caption{Extremal trees $T_2^*(n)$ that achieve the upper bounds.}
\end{figure}

Since $l(T_1^*)\leq l(T_2^*)$, then we have $l(T)\leq l(T_2^*(n)) =n\left \lfloor \frac{(2r-2)n+4}{2r-1}\right\rfloor$ if $n\geq 4r+2$ and $T_2^*(n)$ is the extremal tree that achieves the upper bound. $l(T)\leq l(T_1^*(n)) =\left \lfloor \frac{(2r-2)n+3}{2r-1}\right\rfloor$ if $n\leq 4r+1$  and $T_1^*(n)$ is
 the extremal graph that achieves the upper bound.\qed

Let $T$ be a tree on $n$ vertices. Then $L(T)$ has $n-1$ edges. Since $f(L(T))=l(T)$ and $f(L(T))+\nabla(L(T))=n-1$. Then, we have the following corollary.

\begin{cor}
Let $T$ be a tree on $n$ vertices with diameter $d\geq 4$, then
$$
\left\{\begin{array}{ll}
 \left\lceil\frac{n-d-1}{d-1}\right\rceil \leq \nabla(L(T))\leq n-d-1,&\ for\ d\ is\ even;\\[6pt]
\left\lceil\frac{n-d-2}{d-2}\right \rceil -1 \leq \nabla(L(T))\leq  n-d-1,&\ for\ d\ is\ odd.
\end{array}\right.
$$
\end{cor}

 \begin{thm}
For any connected graph $G$ with $n$ vertices and $m$ edges, if the length of the longest path in $G$ is $p$ and $p\geq 4$, then we have
$$
\left\{\begin{array}{ll}
 m- \left\lfloor \frac{(p-2)n+2}{p-1}\right \rfloor \leq \nabla(L(G))\leq  m-p,& \mbox{ if } p \mbox{ is even;}\\[6pt]
m-\left\lfloor \frac{(p-3)n+4}{p-2}\right \rfloor \leq \nabla(L(G))\leq  m-p, & \mbox{ if } p \mbox{ is odd}.
\end{array}\right.
$$

\end{thm}
\pf  According to Lemma \ref{lem-lf}, we have that $t(L(G))= p(G)=p$ and $f(L(G))=l(G)$. Then it is clear that  $f(L(G))\geq t(L(G))=p$. Since any linear forest in $G$ can be extended to a spanning tree of $G$, then there exists a spanning tree $T$ of $G$ such that $l(G) = l(T)$. Moreover, the maximum linear forest in any spanning tree of $G$ is also a linear forest of $G$. It implies that  for any spanning tree $T$ of $G$, $l(T)\leq l(G)$. Consequently, let $\mathcal{T}$ be the set of all spanning trees of $G$, then $l(G)=\max_{T \in \mathcal{T}} l(T)$. It is easy to see that the diameter of each spanning tree of $G$ is less than or equal to $s$. Moreover, upper bounds on $l(T)$ in Theorem \ref{thm-ex} are all increasing functions of diameters. Thus,
$$
\left\{\begin{array}{ll}
 p \leq l(G)\leq  \left \lfloor \frac{(p-2)n+2}{s-1}\right \rfloor,\ for\ p\ is\ even;\\[6pt]
p \leq l(G)\leq   \left \lfloor \frac{(p-3)n+4}{s-2}\right \rfloor,\ for\ p\ is\ odd.
\end{array}\right.
$$

Thus, the theorem follows.
\qed

 \section{The Decycling Number of Line Graphs of $k$-ary Trees}

 A $k$-ary tree is a rooted tree where within each level every node has either $0$ or $k$ children. The maximum degree of a $k$-ary tree is $k+1$. It follows that line graphs of $k$-ary trees are $K_{k+2}$-free graphs. Since it is often interesting to consider the decycling number of $K_{k+2}$-free graphs. Thus, we consider the decycling number of line graphs of $k$-ary trees in this section.

 Before that, we give a dynamic programming algorithm to find a maximum linear forest in rooted tree. Let $T$ be a rooted tree with root $v^*$. Let $T_1,T_2,...,T_t$ be subtrees of $v^*$ with root $v_1,v_2,\ldots,v_t$, respectively.
 Let $F(T)$ be the edge set of a maximum linear forest in $T$. Let $F'(T)$
 be the edge set of a largest linear forest of $T$ such that the degree of $v^*$ is at most one.

\begin{figure}[ht]
\label{fig-algorithm}
\center
\includegraphics[scale=0.65]{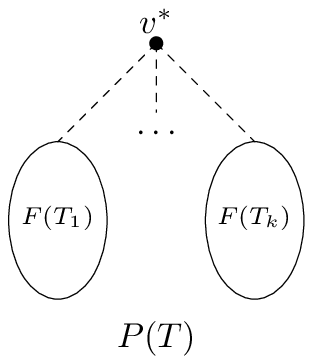}
\hspace{40pt}
\includegraphics[scale=0.65]{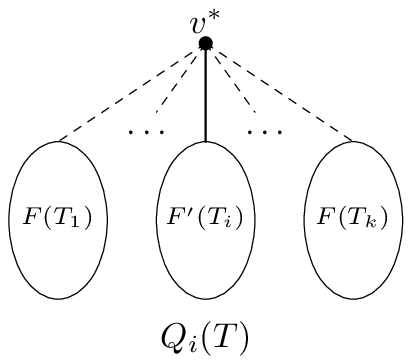}
\hspace{40pt}
\includegraphics[scale=0.65]{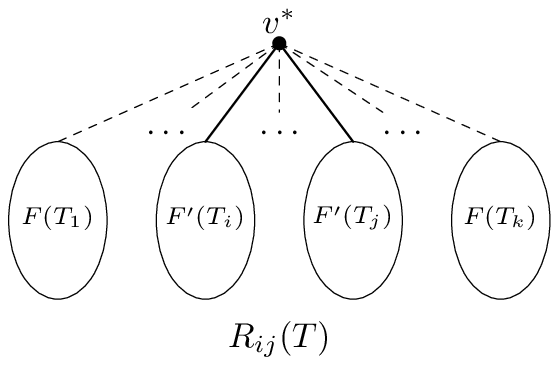}
\caption{The structures of  linear forests $P(T)$, $Q_i(T)$ and $R_{ij}(T)$.}
\end{figure}

Then, we define three kinds of linear forests $P(T)$, $Q_i(T)$ and $R_{ij}(T)$ as follows.
\[ \left\{\begin{array}{l}P(T) =  \cup_{i=1}^k F(T_i),\\[5pt]
Q_i(T) =  \left (\cup_{k\neq i} F(T_k)\right )\cup F'(T_i) \cup \{(v^*,v_i) \},\\[5pt]
R_{ij}(T) =  \left (\cup_{k\neq i,j} F(T_k)\right )\cup F'(T_i)\cup F'(T_j) \cup \{(v^*,v_i),(v^*,v_j)\}.
\end{array}\right .\]
As shown in Fig.6, $P(T)$ is a linear forest in $T$ such that $v^*$ has degree zero;
$Q_i(T)$ is a linear forest such that $v^*$ has degree one and $(v^*,v_i)$ is an edge in the linear forest;
$R_{ij}(T)$ is a linear forest such that $v^*$ has degree two and $(v^*,v_i)$, $(v^*,v_j)$ are  edges in the linear forest.
Let $S(T)$ be the largest linear forest among all $P(T)$, $Q_i(T)$ and $R_{ij}(T)$. Let $S'(T)$ be the largest linear forest among all $P(T)$, $Q_i(T)$.

\begin{lem}\label{lem-algor}
For any tree $T$, $S(T)$ is a maximum linear forest in $T$, $S'(T)$ is a largest linear forest in $T$ such that root $v^*$ has degree at most one.
\end{lem}
\pf
For any tree $T$, if $S(T)$ is not a maximum linear forest in $T$. Let $F(T)$ be a maximum linear forest in $T$.
Let $F[T_t]=F(T) \cap E(T_k)$, for $k=1,2,\ldots,t$.
We can divide the proof into three cases according to the degree of $v^*$ in $F(T)$.

{\bf Case 1}. The degree of $v^*$ in $F(T)$ is zero.
 Clearly,  $F[T_k]$
is a linear forest in subtree $T_k$. Then we have $|F[T_k]|\leq |F(T_k)|$. Therefore, $|F(T)|\leq |P(T)|$.

{\bf Case 2}. The degree of $v^*$ in $F(T)$ is one. Suppose $v^*v_i$ is in $F(T)$. Then  let $F[T_k]=F(T) \cap E(T_k)$, for $r=1,2,\ldots,t$.
Clearly, for $k\neq i$, $F[T_k]$ is a linear forest in subtree $T_k$ and $F[T_k]$ is a linear forest in subtree $T_k$ with degree
of $v^*$ at most one. Thus,  for each $k\neq i$ we have $|F[T_k]|\leq |F(T_k)|$ and $|F[T_i]|\leq |F'(T_i)|$.
Therefore, $|F(T)|\leq |Q_i(T)|$.

{\bf Case 3}. The degree of $v^*$ in $F(T)$ is two. Suppose $v^*v_i$ and $v^*v_j$ are in $F(T)$. Then  let $F[T_k]=F(T) \cap E(T_k)$.
Clearly, for $r\neq i,j$, $F[T_k]$ is a linear forest in subtree $T_k$. $F[T_i]$  is a linear forest in subtree $T_i$ with degree of $v^*$ at most one and $F[T_j]$ is a linear forest in subtree $T_j$ with degree of $v^*$ at most one.
 Thus,  for each $k\neq i,j$ we have $|F[T_r]|\leq |F(T_r)|$. For each $r= i,j$,
 we have $|F[T_i]|\leq |F'(T_i)|$ and $|F[T_j]|\leq |F'(T_j)|$.
Therefore, $|F(T)|\leq |R_{ij}(T)|$.

Combining these cases, we get the conclusion that $F(T)\leq S(T)$, which implies that $S(T)$ is a maximum linear forest in $T$. Similarly, we can prove $S'(T)$ is a largest linear forest in $T$ such that root $v^*$ has degree at most one.
\qed

\begin{thm}\label{lem-perfectkary}
For any $k$-ary  tree $T$ with $n$ vertices, we have
\[
\frac{n+k-1}{k}\leq l(T)\leq \frac{2n-2}{k}.
\]
\end{thm}
\pf
For the lower bound, suppose $T$ has $x$ internal vertices and $y$ leaves, then we have $kx+1=x+y=n$. It follows that $x=\frac{n-1}{k}$ and $y=\frac{(k-1)n+1}{k}$.
By Lemma \ref{lem-lh}, we have $h(T)\leq \frac{(k-1)n+1}{k}-1$. Thus, $l(T)\geq n-h(T)=\frac{n+k-1}{k}$.

For the upper bound, we can divide $n-1$ edges of $T$ into $\frac{n-1}{k}$ groups such that each $k$ edges
 with the same parent are in the same group.
Since the degree of vertices in the linear forest is at most two. Thus, at most two edges in each group are in the linear forest. Therefore,  we get
$l(T)\leq \frac{2n-2}{k}$.\qed

\begin{figure}[ht]
\label{fig-karytree}
\center
\includegraphics[scale=0.85]{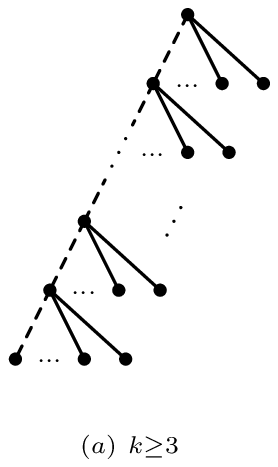}
\hspace{40pt}
\includegraphics[scale=0.85]{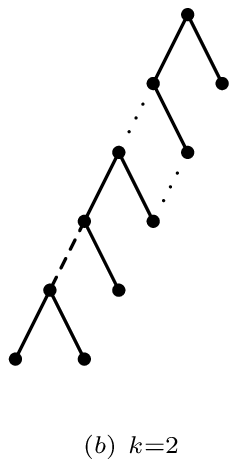}
\caption{Special $k$-ary trees and their maximum linear forests.}
\end{figure}

Let $T$ be a $k$-ary tree on $n$ vertices such that in each layer there is only one node with $k$ children, as shown in Fig. 7. Let $v_i$ be the internal vertex at depth $i$ and root $v_0$ is at depth $0$.
For $k\geq 3$, it is clear that $l(T)= \frac{2n-2}{k}$ as shown in  Fig. 7 (a). For $k=2$, by Theorem \ref{lem-algor} it is easy to check that the linear forest shown in
Fig. 7 (b) is maximum. Thus, $l(T)= \frac{3(n-1)}{4}$ for $\frac{n-1}{2}$ is even; $l(T)= \frac{3n-1}{4}$ for $\frac{n-1}{2}$ is odd.

\begin{cor}
Let $T$ be a $k$-ary tree on $n$ vertices, then
\[
\frac{(k-2)n-k+2}{k}\leq \nabla(L(T)) \leq \frac{(k-1)n-2k+1}{k}.
\]
\end{cor}

A perfect $k$-ary tree is a $k$-ary tree in which all leaves are at the same depth. At last, we obtain the maximum linear forests in perfect $k$-ary trees as follows.

\begin{thm}\label{lem-perfectkary}
For any perfect $k$-ary tree $T$ with $n$ vertices, we have
$$l(T) = \frac{2n-1+(-1)^h}{k+1},$$ where $h=\log_k\left(n(k-1)+1\right)$ is the height of $T$ and leaves are at height 1.
\end{thm}
\begin{figure}[ht]
\label{fig-karytree2}
\center
\includegraphics[scale=0.45]{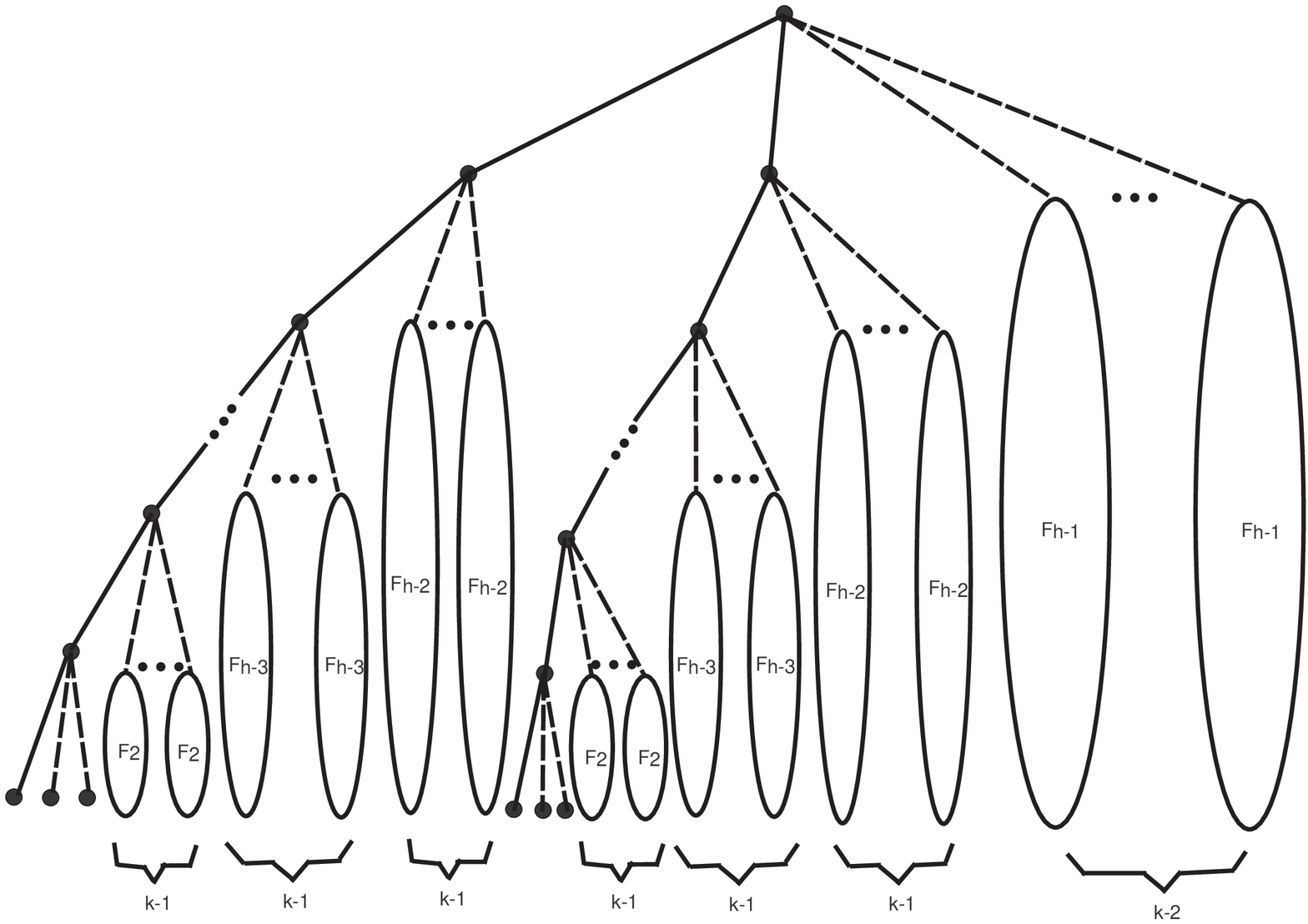}
\caption{A linear forest in perfect $k$-ary trees.}
\end{figure}
\pf Let $h=\log_k\left(n(k-1)+1\right)$ be the height of $T$. We construct a linear forest of $T$ as follows.
Firstly, we choose two vertex-disjoint paths of length $h$ that go from
root to two leaves. Then tree $T$ is decomposed into $k-2$ subtrees of height $h-1$,  $2(k-1)$ subtrees of height $h-2$,
$2(k-1)$ subtrees of height $h-3$, \ldots, $2(k-1)$ subtrees of height $3$ and $2(k-1)$ subtrees of height $2$ as shown in Fig. 8.
Then for each subtree, we  choose two vertex-disjoint paths that go from the root to the  two leaves again. Do it recursively, then  a linear forest of $T$ is created.

Let $F_h$ be the edge set of the obtained linear forest in $T$ with height $h$ and
let $f_h$ be the cardinality of $F_h$.
Then it is easy to see that $f_1=0$ and $f_2=2$. According to the recursive construction of the obtained linear forest, we have
\[
f_h=(k-2)f_{h-1}+2(k-1)\sum_{i=1}^{h-2} f_i+2(h-1),\]
and
\[
f_{h-1}= (k-2)f_{h-2}+2(k-1)\sum_{i=1}^{h-3} f_i+2(h-2).
\]
Combining the two equations, we get a recursive relation as follows.
\[f_h=(k-1)f_{h-1}+kf_{h-2}+2.\]

By the technique of generating functions, we can derive a formula for $f_h$ as follows.
\[
f_h=\frac{2n-1+(-1)^h}{k+1}.
\]
Thus,
$l(T)\geq f_h = \frac{2n-1+(-1)^h}{k+1}$.

Now we prove that $F_h$ is a maximum linear forest in $T$. Let $v^*$ be the root of $T$ and $T_1$, $T_2$, \ldots, $T_t$ be
the $t$ subtrees of $v^*$ with root $v_1,v_2,\ldots,v_t$. Since $T$ is a perfect $k$-ary tree, each subtree $T_k$ is identical to a perfect $k$-ary tree of height $h-1$. Suppose $v^*v_1$ and $v^*v_2$ are in $F_h$.  Let $F'_{h-1}$ be the subset of $F_h$ in subtree $T_1$ and let $f'_{h-1}$ be the cardinality of $F'_{h-1}$.
Then $F'_{h-1}$ is a linear forest in $T_1$ such that $deg(v_1)=1$.

We claim that in a perfect $k$-ary tree of height $h$, $F_h$ is a maximum linear forest and $F'_h$ is a largest linear forest such that the degree of the root is at most one. We prove the claim by induction on $h$. For $h=1$ and $h=2$, it is easy to check $F_1$, $F_2$ are maximum linear forests and $F'_1$, $F'_2$ are largest linear forests with degree of root at most one, where $F_1,F'_1$ are empty sets.
Suppose the claim is true for perfect $k$-ary tree with height $h-1$. Let $T$ be a perfect $k$-ary tree with height $h$. Define
\[ \left\{\begin{array}{ll}P(T) &=  \cup_{i=1}^k F(T_i),\\[5pt]
Q_i(T) &=  \left (\cup_{k\neq i} F(T_k)\right )\cup F'(T_i) \cup \{(v^*,v_i) \},\\[5pt]
R_{ij}(T) &=  \left (\cup_{k\neq i,j} F(T_k)\right )\cup F'(T_i)\cup F'(T_j) \cup \{(v^*,v_i),(v^*,v_j)\}.
\end{array}\right .\]

By induction hypothesis, each $F(T_i)$ is identical to $F_{h-1}$ and each $F'(T_i)$ is identical to $F'_{h-1}$.
Then
\[ \left\{\begin{array}{ll}|P(T)|&=kf_{h-1},\\[5pt]
|Q_i(T)|&=(k-1)f_{h-1}+f'_{h-1}+1,\\[5pt]
|R_{ij}(T)|&=(k-2)f_{h-1}+2f'_{h-1}+2.
\end{array}\right .\]
It is easy to see that each $T_i$ has $\frac{n-1}{k}$ vertices. Since $F_h$ is also consist of $k-2$  $F_{h-1}$'s, two $F'_{h-1}$'s and two extra edges, then we have $|R_{ij}(T)|=f_h$. It follows that $f'_{h-1}=\frac{1}{2}\left(f_h-(k-2)f_{h-1}-2\right)$. Similarity, we have $|Q_i(T)|=f_h'$.
Then,
\begin{align*}
f'_{h-1}+1-f_{h-1}&=\frac{1}{2}\left(f_h-(k-2)f_{h-1}-2\right)+1-f_{h-1}\\[5pt]
&=\frac{1}{2}(f_h-kf_{h-1})\\[5pt]
&=\frac{1}{2}\left(\frac{2n-1+(-1)^h}{k+1}-k\frac{2\frac{n-1}{k}-1+(-1)^{h-1}}{k+1}\right)\\[5pt]
&=1-(-1)^{h-1}\geq 0.\\[5pt]
\end{align*}
It follows that  $|P(T)|\leq |Q_i(T)|\leq |R_{ij}(T)|$. By Theorem \ref{lem-algor}, we see that $R_{ij}(T)$ is a maximum linear forest in $T$ and $Q_i(T)$ is a largest linear forest such that degree of the root is at most one, which are exactly $F_h$ and $F_h'$. Therefore, we prove the claim and $F_h$ is a maximum linear forest of $T$. Thus, we conclude that $l(T) = f_h= \frac{2n-1+(-1)^h}{k+1}$.
\qed

\begin{cor}
For any perfect $k$-ary tree $T$ with $n$ vertices, the decycling number of $L(T)$ is
\[
\nabla(L(T))=\frac{2n-1+(-1)^h}{k+1},
\]
where $h=\log_k\left(n(k-1)+1\right)$ is the height of $T$.
\end{cor}

{\noindent \bf Acknowledgements.}

The work was supported by National Natural Science Foundation of China (No.11671299, No.61502330, No.61472465, No.61562066).

\end{document}